\documentstyle[12pt]{article}
  \newcommand{\const}{\rm const}
  
  \newcommand{\supp}{\rm supp}

 \newcommand{\vraisup}{\rm vraisup}
  \newcommand{\Dom}{\rm  Dom}

\begin{document}

   \begin{center}

  \   {\bf   Equivalence between tails, Grand Lebesgue Spaces and Orlicz }\\

\vspace{4mm}

   \  {\bf  norms for random variables  without Cramer's condition  }\\

\vspace{4mm}

  \   {\bf Kozachenko Yu.V.,  Ostrovsky E., Sirota L.}\\

\vspace{4mm}

Department of Probability Theory and Statistics, Kiev State University, Kiev, Ukraine,\\
e-mails: yvk@univ.kiev.ua \\
 ykoz@ukr.net

\vspace{4mm}

 Israel,  Bar-Ilan University, department  of Mathematic and Statistics, 59200, \\

\vspace{4mm}

e-mails: eugostrovsky@list.ru \\
sirota3@bezeqint.net \\

\vspace{5mm}

  {\bf Abstract} \\

\vspace{4mm}

 \end{center}

  We offer in this paper the non-asymptotical pairwise bilateral exact up to multiplicative constants interrelations between the
 tail behavior, moments (Grand Lebesgue Spaces) norm and Orlicz's norm for random variables  (r.v.),
which does not satisfy in general case the  Cramer's condition.\par

\vspace{4mm}

  \par

\vspace{4mm}


\vspace{5mm}

\section{ Definitions.  Notations. Previous results.  Statement of problem.}

 \vspace{4mm}

 \ Let  $ \   (\Omega, \ B, \ {\bf P} )  \ $ be certain probability space with non - trivial probability measure $ \ {\bf P} \ $ and
correspondent expectation $ \ {\bf E}, \ $  and let $  \ \xi = \xi(\omega), \ \omega \in \Omega \ $ be numerical
valued random variable.  We denote as usually by $ \ |\xi|_p, \ p \in [1, \infty] \ $  its classical Lebesgue - Riesz
$ \ L_p = L(p)  = L_p(\Omega) $ norm

$$
|\xi|_p := \left[ {\bf E} |\xi|^p  \right]^{1/p}, \ 1 \le p < \infty;  \
|\xi|_{\infty} := \vraisup_{\omega \in \Omega} |\xi(\omega)|.
$$
 \ The so - called {\it  tail - function } $ \   T_{\xi}(u), \ u \ge 0  \ $  for this random variable $ \ \xi \ $ is defined by a formula

$$
T_{\xi}(u) \stackrel{def}{=} \max \left\{  \ {\bf P} (\xi > u), \ {\bf P}(\xi <-u) \   \right\}, \ u \ge 0.\eqno(1.1)
$$

 \ An equivalent version:

$$
\overline{T}_{\xi}(u) := {\bf P} (|\xi| > u), \ u \ge 0. \eqno(1.1a)
$$

 \ Obviously,

$$
 \overline{T}_{\xi}(u)   \le T_{\xi}(u) \le  2 \ \overline{T}_{\xi}(u), \ u \ge 0,
$$
 the equivalence. \par

\vspace{4mm}

 {\bf The aim of this report is to establish the reciprocal non-asymptotic interrelations separately mutually
possibly exact up to multiplicative constant between tail functions, suitable Orlicz and Grand Lebesgue Spaces norms
for random variables. }\par

 {\bf   We do not suppose that the considered in this article r.v. satisfy the famous Cramer's condition:  }

$$
\exists \epsilon_0 > 0 \  \forall \lambda: \ |\lambda| < \epsilon_0 \  \Rightarrow  {\bf E} \exp(\lambda \xi) < \infty,
$$

 \ {\bf   in contradiction  with previous works, see, for example, works [2], [6], [7], [8], [19]}. \par

\vspace{4mm}

 {\it Throughout this paper, the letters  $ \ C, C_j(\cdot) \ $ etc. will denote a various positive finite
constants which may differ from one formula to the next even within a single string
of estimates and which does not depend on the essentially variables  $  \ p, x, \lambda, y, u \ $ etc. \par

 \ We make no attempt to obtain the best values for these constants.}\par

 \ The immediate predecessor of offered report is the article [7],  in which was considered the case when
the considered  r.v.  satisfy the famous Cramer's condition.  See also  [6], [8], chapters 1,2  etc. \par

 \  We will use in this report  in general at the same techniques  as in [7]. \par

 \ Recall that the  so-called Young-Fenchel transform $ \  g \to g^* \  $ of arbitrary real valued function $  \  g = g(x) \ $
is defined as follows

$$
g^*(y) \stackrel{def}{=} \sup_{x \in \Dom(g)} (x y - g(x)). \eqno(1.2)
$$
 \ The  symbol $ \ \Dom(g) \ $ denotes as ordinary the domain of definition (in particular, finiteness) of the function $ \ g(\cdot). $\par

 \ Let us bring some used further examples. Define the function

$$
\phi_{m,L} = \phi_{m,L}(\lambda) :=  m^{-1} \ \lambda^m \ L(\lambda), \ \lambda > 0, \ m = \const > 1, \eqno(1.3)
$$
where $  \ L = L(\lambda) $ is positive {\it  slowly varying }  at infinity, i.e. as $ \ \lambda \to \infty \ $
function. Then as $  \ x \to \infty $

$$
\phi^*_{m,L}(x) \sim  (m')^{-1} \ x^{m'} \ L^{-1/(m-1)} \left(  x^{1/(m-1)}   \right), \eqno(1.4)
$$
and as ordinary  for arbitrary value $ \ m > 1 $
$$
  m' \stackrel{def}{=} \frac{m}{m-1}.
$$

  \ If   for instance  $ \ L(\lambda) = [\ln (\lambda + e)] ^r, \ r = \const \in R, \ $  i.e.

$$
\phi(\lambda) := \phi_{m,r}(\lambda) = m^{-1} \ \lambda^m \ [\ln (\lambda + e)]^r, \ \lambda \ge 0, m = \const > 1, \ r \in R,
$$
then as $ \ x \to \infty $

$$
\phi_{m,r}^*(x) \sim (m')^{-1} \ x^{m'} \ [\ln ( x + e)]^{-r/(m-1)}, \eqno(1.5)
$$
 see, e.g.  [18],  p. 40 - 42. \par
 \ Analogously if  $ \ \phi(\lambda) := \phi_{m,r,q}(\lambda) = $

$$
 m^{-1} \ \lambda^m \ [\ln ( \lambda + e)]^r \  [\ln \ln (\lambda + e^e)]^q,
 \ \lambda > 0, m = \const > 1, \  r,q \in R,
$$
then similarly  as $ \ x \to \infty $

$$
\phi_{m,r,q}^*(x) \sim (m')^{-1} \ x^{m'} \ [\ln (x+e)]^{-r/(m-1)} \ [\ln \ln( x + e^e)]^{-q/(m-1)}. \eqno(1.6)
$$

  \ More generally, if

$$
L(\lambda) = [\ln \lambda]^r \ M(\ln \lambda),
$$
where  $  \ M = M(\lambda) \ $ is positive   slowly varying  function  as $  \ \lambda \to \infty,  $ then  as $ \ x \to \infty $

$$
\phi^*_{m, L}(x) \sim  (m')^{-1} \ x^{m'} \ [\ln x]^{ - r/(m-1)  } \ M^{-1/(m-1)} (\ln x). \eqno(1.7)
$$

\vspace{4mm}

 \ The case $  \ m = 1 \ $ is more complicated. Define the $  \  \Psi \  $ function $ \ \psi^{(L)} = \psi^{(L)}(p) \   $ as follows

$$
\psi^{(L)}(p)  \stackrel{def}{=} \frac{p}{ L(p)}, \eqno(1.8)
$$
where as before  $  \ L = L(\lambda) \ $ is positive  continuous  slowly varying  function  as $  \ \lambda \to \infty $ tending to
infinity as $ \ \lambda \to \infty. \ $   Let also the r.v. $ \ \xi \ $ be from the Grand Lebesgue Space $ \ G\psi^{(L)} \  $
with unit norm:

$$
||\xi||G\psi^{(L)}  \stackrel{def}{=} \sup_{p \ge 1} \left\{ \frac{|\xi|_p}{ \psi^{(L)}(p)}  \right\}   = 1,
$$
then by the direct definition of these norms

$$
|\xi|_p \le \frac{p}{L(p)}, \ p \ge 1. \eqno(1.9)
$$

 \ We deduce by means of Tchebychev-Markov inequality

$$
T_{\xi}(x) \le \exp \left(  - C(L) \ x \ \ln L(x)  \right). \eqno(1.10)
$$

 \ Conversely, let the estimate (1.10) be a given for some r.v. $ \ \xi \ $ and some such a  positive function $  \ L = L(\cdot) $ for which
$    \   L(x) \uparrow \infty,  \  $ then

$$
|\xi|_p \le C(L) \  \frac{p}{L(p)}, \ p \ge 1. \eqno(1.11)
$$

\vspace{4mm}

\section{Grand Lebesgue Spaces (GLS). }

 \vspace{4mm}

 \ Let   $ \   (\Omega, B, {\bf P})  \ $ be again at the same  probability space.
 \ Let also  $  \psi = \psi(p), \ p \in [1, b), \ b = \const \in (1,\infty]  $ (or   $ p \in [1,b] $ ) be certain bounded
from below:  $  \ \inf \psi(p)  > 0 $ continuous inside the  {\it  semi - open} interval $   \ p \in [1, b) $ numerical function
such that the function

$$
h(p) = h[\psi](p) \stackrel{def}{=} p \ \ln \psi(p) \eqno(2.0)
$$
is convex.\par

 \ An important example. Let  $ \ \eta \ $ be a random variable such that there exists $ \  b = \const > 1   \ $  so that
$ \  |\xi|_b < \infty.  \ $  The {\it  natural } $ \ G\Psi \ $  function $ \  \psi_{\eta} = \psi_{\eta}(p)  \  $  for the r.v. $ \ \eta \ $
is defined by a formula

$$
\psi_{\eta}(p) \stackrel{def}{=} |\eta|_p.
$$

 \ We can and will suppose
$   \ b = \sup \{p, \psi(p) < \infty\},  \ $ so that  $ \  \supp \ \psi = [1, b) \  $  or $ \ \supp \ \psi = [1, b]. \ $ The set of all such a
functions will be denoted by  $ \ \Psi(b) = \{  \psi(\cdot)  \}; \ \Psi := \Psi(\infty).  $\par

\vspace{4mm}

 \  {\it  We will consider in this article only the case when $ \ b = \infty; $ i.e. $ \   \Psi := \Psi(\infty). \  $ }\par

\vspace{4mm}

 \ By definition, the (Banach) Grand Lebesgue Space  \ (GLS)  \ space   $  \ G\psi = G\psi(b)  $ consists on all the
numerical  valued  random variables
(measurable functions) $ \zeta  $ defined on our measurable space  and having a finite norm

$$
||\zeta|| = ||\zeta||G\psi \stackrel{def}{=} \sup_{p \in [1,b)} \left\{ \frac{|\zeta|_p}{\psi(p)} \right\}. \eqno(2.1)
$$

\  The function $ \  \psi =\psi(p) \  $ is named {\it  generating function } for the Grand Lebesgue Spaces. \par

 \ These spaces  are Banach functional space, are complete, and rearrangement invariant in the classical sense, see [1],  chapters 1, 2;
and were  investigated in particular in  many  works, see e.g.  [3], [4], [5], [6], [7], [8], [15], [16], [17].
 We refer here  some  used in the sequel facts about these spaces and supplement more. \par

\vspace{4mm}

\  It is known that if  $  \ \zeta \ne 0,  $ and $ \  \zeta \in G\psi(b), \ $ then

  $$
  T_{\zeta} ( y) \le \exp \left( \ - h_{\psi}^* (\ln ( y/||\zeta||) ) \ \right), \ y \ge ||\zeta||, \eqno(2.2)
  $$
where

$$
h(p) = h[\psi](p)  \stackrel{def}{=} p \ \ln \psi(p), \  1 \le p < b.
$$

 \ Namely, let $ \  ||\zeta||_{G\psi(b) } = 1;  \ $  therefore by means of Tchebychev-Markov inequality

$$
 T_{\zeta} ( y)  \le  \frac{\psi^p(p)}{y^p} =\exp \left(  - p \ln y + p  \psi(p)    \right),
$$
following

$$
 T_{\zeta} ( y)  \le    \inf_{ p \in [1,b) } \exp \left(  - p \ln y + p \ \psi(p)    \right) =
\exp \left( \ - h[\psi]^* (\ln (y/||\zeta||) ) \ \right), \ y \ge e \cdot ||\zeta||.
$$

  \ Conversely, the last inequality may be reversed in the following version: if the r.v. $  \  \zeta \  $  {\it satisfies the Cramer's
condition} and

$$
{\bf P}(|\zeta| > y) \le \exp \left(-h_{\psi}^* (\ln (y/K) \right). \ y \ge e \cdot K, \ K = \const \in (0,\infty),
$$
and if the function $ h_{\psi}(p), \ 1 \le p < \infty  \ $  is positive, continuous, convex and such that

$$
\lim_{p \to \infty}  \psi(p)/p = 0,
$$
then  $ \zeta \in G\psi,  $ herewith $ \ ||\zeta||  \le C(\psi) \cdot K  $  and  conversely

$$
||\zeta||G\psi \le C(\psi) K  \le C_2(\psi) |\zeta||G\psi, \ 0 < C_1(\psi) < C_2(\psi) < \infty. \eqno(2.3)
$$

 \ Introduce the  following {\it exponential} Young-Orlicz function

$$
N_{\psi}(u) = \exp \left(h_{\psi}^* (\ln |u|) \right),  \ |u| \ge 1; \ N_{\psi}(u) = C u^2, \ |u| < 1,
$$
and the correspondent Orlicz norm will be denoted by $ \ ||\cdot||L \left(N_{\psi} \right) =  ||\cdot||L (N). \ $ It was  done

$$
||\zeta||G\psi \le C_1 ||\zeta||L(N)  \le C_2 ||\zeta||G\psi, \ 0 < C_1 < C_2 < \infty. \eqno(2.4)
$$

 \ If for instance $  \ \psi(p) = \psi_m(p)\stackrel{def}{=} p^{1/m}, \ p \in [1, \infty), $ where $  \ m = \const > 1, \ $  then

$$
 0 \ne \xi \in G\psi_m \Leftrightarrow \ T_{\xi}(u) \le \exp \left(-C(m) u^m \right).
$$
  \ Define also the correspondent Young-Orlicz function

$$
N_m(u) := \exp \left( |u|^m \right), \ |u| \ge 1; \ N_m(u) = C u^2, \ |u| \le 1.
$$
 \  The relation (2.3) means in addition in this case

$$
||\zeta||G\psi_m \le C_1(m) ||\zeta||L(N_m)  \le C_2 ||\zeta||G\psi_m, \ 0 < C_1(m) < C_2(m) < \infty.  \eqno(2.5)
$$

  \ Notice that  in the case when $  \ m \in (0,1) \ $ the correspondent random variable $ \ \xi \ $ does not satisfy the Cramer's
condition. We intend to generalize the last propositions further on the case just in particular $  \ m \in (0,1]. $

 \ Define as an example  the following {\it  degenerate }  $ \  G\Psi \  $ function

$$
\psi_{(r)}(p) = 1, \ 1 \le p \le r;  \ \psi_{(r)}(p) = \infty, \ p > r; \  r = \const > 1.
$$

 \ The $ \ G\psi_{(r)}   \ $ norm of an arbitrary  r.v. $ \   \eta \  $ is quite equivalent to the classical Lebesgue-Riesz
$ \ L_r \ $ norm

$$
||\eta|| G\psi_{(r)} = |\eta|_r. \eqno(2.6)
$$
 \ Thus, the Grand Lebesgue Spaces are direct generalizations of the Lebesgue-Riesz spaces. \par

\vspace{4mm}

\section{Auxiliary estimates from the  saddle-point  method.}

 \vspace{4mm}

 \ We must investigate  in advance one interest and needed further integrals. Namely, let $  \  (X, M, \mu), \ X \subset R  \  $
 be non-trivial measurable space with non-trivial sigma finite measure  $ \ \mu. \ $ \par

 \ We assume at once  $ \ \mu(X) = \infty, \ $ as long as the opposite case is trivial for us.
We intend to estimate for sufficiently greatest values of real parameter $ \ \lambda \ $,  say $  \ \lambda > e, \ $  the
following integral

$$
I(\lambda) := \int_X e^{  \lambda x -  \zeta(x)  } \ \mu(dx).  \eqno(3.1)
$$
assuming of course its convergence for all the sufficiently great values of the parameter $ \  \lambda. \ $  The offered below
estimates may be considered as a some generalizations of the  saddle-point  method. \par

 \ Here  $   \ \zeta = \zeta(x)   \  $ is non-negative measurable function, not necessary to be convex. \par

 \ We  represent now two methods for {\it upper} estimate $  I(\lambda) $ for sufficiently greatest values of the real parameter $ \ \lambda. \ $ \par

  \ Note first of all that if in contradiction the measure   $  \ \mu \ $ is finite: $ \  \mu(X) = M \in (0, \infty); \ $ then the
integral  $ \ I(\lambda) \ $ allows a  very simple estimate

$$
 I(\lambda)  \le M \cdot \sup_{x \in X}  \exp \left( \lambda x - \zeta(x)    \right) =
M \cdot \exp \left( \zeta^*(\lambda)  \right).  \eqno(3.2)
$$

 \ Let now $   \ \mu(X) = \infty $ and let $ \   \epsilon = \const \in (0,1); \ $
 let us introduce the following  auxiliary integral

$$
K(\epsilon) := \int_X e^{-  \epsilon \zeta(x) } \mu(dx). \eqno(3.3)
$$

 \ It will be presumed its finiteness at last for some positive value  $ \   \epsilon_0 \in (0,1); \ $
 then   $ \  \forall \epsilon \ge \epsilon_0 \ \Rightarrow K(\epsilon) < \infty. \  $

Then the following measures are probabilistic:

$$
\nu_{\epsilon}(A ) := \frac{\int_A \exp( - \epsilon \zeta(x)) \ \mu(dx)}{K(\epsilon)},   \ \epsilon \ge \epsilon_0. \eqno(3.4)
$$

 \ We have

$$
\frac{I(\lambda)}{K(\epsilon)} = \int_X \exp(\lambda x - (1 - \epsilon) \zeta(x)) \  \nu_{\epsilon}(dx) \le
$$

$$
\exp \{ \sup_{x \in X} [  \lambda x - (1 - \epsilon) \zeta(x)   ]    \} =
\exp \left\{ (1 - \epsilon) \zeta^* \left(  \frac{\lambda}{1 - \epsilon} \right)  \right\}.
$$

 \ Following,

$$
I(\lambda) \le K(\epsilon) \cdot  \exp \left\{ (1 - \epsilon) \zeta^* \left(  \frac{\lambda}{1 - \epsilon} \right)  \right\}
\eqno(3.5)
$$
and hence:

\vspace{4mm}

{\bf  Theorem 3.1 }  We assert actually under formulated here conditions, in particular, the condition of the finiteness
of $  \ K(\epsilon) $ for some value $  \ \epsilon_0 \in (0,1): $

$$
I(\lambda) \le \inf_{\epsilon \in (0,1)}
\left[ K(\epsilon) \cdot  \exp \left\{ (1 - \epsilon) \zeta^* \left(  \frac{\lambda}{1 - \epsilon} \right)  \right\} \right].
\eqno(3.6)
$$

 \vspace{4mm}

 \ We can detail the  choice  of the value  $  \  \epsilon \    $ in the estimates  (3.5) - (3.6). Namely,  denote

$$
 \theta = \theta( \lambda ) := \frac{c_1}{ \lambda \ \zeta^{*'} (2\lambda)},  \ \lambda \ge \lambda_0 = \const > 0.  \eqno(3.7)
$$

  \ The value $ \ \lambda_0 \ $ is selected such that $  \    \theta(\lambda) \le 1/2, \ \lambda \ge \lambda_0.   \ $ Then

$$
\frac{\lambda}{1 - \epsilon} \le \lambda(1 + 2 \epsilon),
$$
and we have  taking into account the convexity of the function $ \  \zeta^*(\cdot)  \  $ and denoting
$  \   \phi(\lambda) = \zeta^*(\lambda): \  $

$$
\phi \left( \frac{\lambda}{1 - \theta}  \right) \le \phi(\lambda + 2 \lambda \theta)  \le
$$

$$
\phi(\lambda) +2 \theta \lambda \ \phi'(2 \lambda) \le  c_2 + \phi(\lambda).
$$

 \ To summarize:

\vspace{4mm}

$$
I(\lambda) \le c_2 \  K(\theta(\lambda))  \  \exp(\zeta^*(\lambda)). \eqno(3.8)
$$

 \ As regards  the function $  \ K = K(\theta(\lambda)), \ $   note that if  $ \  X   = R^+, \ \mu(dx) = dx, \ $ and if

$$
\zeta(x) \ge c_4 \ x, \ x \ge 0, \eqno(3.9)
$$
then

$$
K(\theta(\lambda)) \le c_5 \ \lambda \ \zeta^{*'}(2 \lambda),
$$
hence

$$
I(\lambda) \le c_6 \  \lambda \ \zeta^{*'}(2 \lambda) \cdot  \exp(\zeta^*(\lambda)), \ \lambda > \lambda_0. \eqno(3.10)
$$

 \ If in turn instead (3.9) there holds

$$
\zeta(x) \ge c_7 \ x^{\alpha}, \ \alpha = \const > 0,  \ X = R_+, \ \mu(dx) = dx,
$$
then

$$
I(\lambda) \le c_8 \  \left[ \lambda \ \zeta^{*'}(2 \lambda) \right]^{1/\alpha} \cdot  \exp(\zeta^*(\lambda)), \
\lambda > \lambda_0. \eqno(3.11)
$$

\vspace{4mm}

{\bf Theorem 3.2.}  Suppose in addition  $   \ X = (a, \infty), \ a = \const \in R,\ $  or $ \  X = R, \ $
and that

$$
\exists C = \const \in (0, \infty), \exists \alpha = \const > 1 \ \Rightarrow
\zeta(x) \ge C x^{\alpha}, \ x \ge 1. \eqno(3.12)
$$
 \ Then  there exists a finite positive constant $ \  C = C(\zeta, a)  \  $ such that for sufficiently values $  \lambda, $ say  for
$  \ \lambda \ge 1 $

$$
I(\lambda) \le \exp \left(  \zeta^*(C \lambda)  \right).   \eqno(3.13)
$$

 \vspace{4mm}

  \ {\bf Proof,}   in particular, the finiteness of $ \ K(\epsilon), \ \epsilon \in (0,1) \ $ contains in fact in [9],
chapter 2.1. \par

\vspace{4mm}

 \  We represent now an opposite method, which was introduced in particular case in  [7], [8], sections 1.2.
 Indeed,  let $ \   \gamma = \const \in (0, 1).  \  $ We apply  the Young's inequality

$$
\lambda x \le \zeta(\gamma x) + \zeta^*(\lambda/\gamma),
$$
therefore

$$
I(\lambda) \le e^{  \zeta^*(\lambda/\gamma)  } \cdot \int_X e^{\zeta(\gamma x) - \zeta(x)    } \ \mu(dx) =
R(\gamma)  \ e^{  \zeta^*(\lambda/\gamma)  },  \eqno(3.14)
$$
where

$$
R(\gamma) := \int_X e^{\zeta(\gamma x) - \zeta(x)    } \ \mu(dx).  \eqno(3.15)
$$

  \ We obtained  really the following second estimate. \par

\vspace{4mm}

{\bf Lemma 3.2.}

$$
I(\lambda) \le \inf_{\gamma \in (0,1)} \left[ R(\gamma)  \ e^{  \zeta^*(\lambda/\gamma)  }   \right].  \eqno(3.16)
$$

\vspace{4mm}

\section{Main results: connection between tail behavior and Grand Lebesgue Space norm.}

 \vspace{4mm}

  \ {\it  Statement of problem:}  given a tail function $ \ T_{\xi} (y) \ $ for the certain (non-zero) random variable
$ \  \xi \ $  of the form

$$
T_{\xi} (y) \le \exp \left(  - h^*[\psi](\ln y)   \right),  \ y \ge 1, \eqno(4.1)
$$
 where $  \ \psi(\cdot) \in G\Psi. $  It is  required to prove $ \  \xi \in G\psi, \  $ or on the other words to obtain an estimate
of the form $ \ ||\xi ||G\psi < \infty.\ $ \par

 \ Recall that the inverse conclusion: $  \ ||\xi||G\psi = 1 \ \Rightarrow \ $  (4.1) is known, see (2.2).  \par

 \ So, let the estimate (4.1) be a given. We have  for the values $  \  p \ge e $

$$
p^{-1} |\xi|_p^p \le \int_0^{\infty} x^{p-1}  \exp \left(  - h^*[\psi](\ln x)   \right) \ dx =
$$

$$
\int_{-\infty}^{\infty} \exp(p \ y  - h^*(y) ) \ dy. \eqno(4.2)
$$

 \   It remains to use the proposition of theorem 3.1. \par

\vspace{4mm}

{\bf Theorem 4.1.}   Suppose

$$
 C(h) := \sup_{p \in [1, \infty)} \left[ \ h^{* `} [\psi](p) \   \right]^{1/p} < \infty. \eqno(4.3)
$$
 \ If the r.v. $ \ \xi \ $ satisfies the inequalities (4.1) and (4.3), then  $ \  \xi \in G\psi: \  $

$$
 ||\xi ||G\psi \le 2 \ C[h] \ e^{1/e} < \infty. \eqno(4.4)
$$

 \ {\bf Proof.}  It is sufficient  to note that the function $ \  p \to h[\psi(p)]  \  $ is continuous and convex and that

$$
\left(h^* \right)^* = h^{**} = h
$$
by virtue of theorem of Fenchel-Moreau. \par

\vspace{4mm}

\  Let us bring some examples.  \\

\vspace{4mm}

{\bf Example 4.1.}  Put as before

$$
\psi_m(p) = p^{1/m},
$$
but here  $  \ m = \const \in (0, \infty).  $ Let $ \ \xi \in G\psi_m \ $ and $ \ ||\xi||G\psi_m = 1.  \ $ \par
 \ Note that in the case $ \   m \in (0,1) \ $ the r.v. $ \  \xi \  $ does not satisfy in general case the Cramer's condition.  But we
conclude on the basis of theorem 3.1  $ \ ||\xi||G\psi_m  \in (0,\infty) \ \Longleftrightarrow $

$$
  \exists \ C(m) \in (0,\infty),  \ T_{\xi}(u) \le \exp \left(  - C(m) \ u^m \right), \ u \ge 0. \eqno(4.5)
$$

 \ More precisely, if  $ \ ||\xi||G\psi_m = 1,  \ $ then

$$
T_{\xi}(u) \le \exp \left( - (me)^{-1} \ y^m \right), \ y > 0.
$$

 \ Inversely, assume

$$
T_{\xi}(u) \le \exp \left( - \ y^m \ \right), \ y > 0.
$$
 \ Then it follows from theorem 3.1

$$
||\xi||G\psi_m \le e^{ m + 1/e }
$$
or equally

$$
|\xi|_p \le e^{ m + 1/e } \ p^{1/m}, \ p \ge 1.
$$

\vspace{4mm}

\ Let us consider a more general case, indeed, introduce as above the following $ \  \Psi \ $ function

$$
\psi_{m,L}(p) \stackrel{def}{=} p^{1/m} \ L(p), \ m = \const > 0, \eqno(4.6)
$$
where $ \  L = L(p), \ p \ge 1 \  $ is some positive continuous slowly varying as $ \ p \to \infty \ $ function. We impose
for simplicity the following condition on this function:

$$
\forall \theta > 0 \ \Rightarrow \sup_{p \ge 1} \left[ \frac{L(p^{\theta})}{L(p)} \right] =: C(\theta) < \infty. \eqno(4.7)
$$

 \ This condition is satisfied, if for example $ \  L(p) = [ \ln (p+1)]^r, \ r = \const.  \  $\par
 \ It follows again from theorem 3.1 that the r.v. $ \  \xi \ $ belongs to the space $ \  G\psi_{m,L}:  \  $

$$
||\xi||  G\psi_{m,L} = \sup_{p \ge 1} \left[\frac{|\xi|_p}{\psi_{m,L}(p)}\right]  = 1  \eqno(4.8)
$$
if and only if

$$
T_{\xi}(y) \le \exp \left( - C(m,L) \ y^m / L(y)   \right), \ y \ge e. \eqno(4.9)
$$

\vspace{4mm}

\ As a particular case: define the $ \  \Psi \ -  \ $ function

$$
\psi_{(m,r)}(p) := p^{1/m} \ \ln^{r}(p + 1), \ p  \ge 1; \ m = \const > 0, \ r = \const \in R. \eqno(4.10)
$$

 \ The random variable $ \  \xi \ $ belongs to the space $ \  G\psi_{(m,r) }: \ $

$$
||\xi||G\psi_{(m,r)} =  \sup_{p \ge 1} \left[ \ \frac{|\xi|_p}{\psi_{(m,r)}(p)} \ \right]  = l \in (0, \infty)  \eqno(4.11a)
$$
if and only if

$$
T_{\xi}(u) \le \exp \left( - C(m,r) \ (u/l)^m \ \ln^{ -r  }(u/l) \right), \ u \ge e \ l. \eqno(4.11b)
$$

\vspace{4mm}

 \ \ \ {\bf Example 4.2.} A boundary case. \par

\vspace{4mm}

 \  We introduce the following $ \ G\Psi \ $ function

$$
\psi^{(s)}(p) = p \ (\ln(p+1))^s, \ s = \const \in R, \ p \in [1, \infty).  \eqno(4.12a)
$$
 \ Then  the non - zero r.v. $  \  \nu \ $ belongs to the $ \ G \psi^{(s)} \ $ space if and only if

$$
T_{\nu}(y) \le \exp \left( - C(s) \ y \ \ln^{-s}(y+1)  \right), \ y \ge 0.  \eqno(4.12b)
$$

 \ Note that the r.v.  $ \ \nu \ $ satisfies the Cramer's condition if and only if $ \ s \le 0. \ $  The case $  \ s = 0 \ $
correspondent to the exponential distribution for the r.v. $ \ \nu; \ $  the case $ \ s = - 1 \ $ take place in particular when  the r.v.
$ \ \nu \ $ has a Poisson distribution, which obey's  but the exponential moments. \par

\vspace{4mm}

{\bf Example 4.3.}

\vspace{4mm}

 \  Let us consider the following $  \   \psi_{\beta}(p) \  $ function

$$
  \psi_{\beta,C}(p)  :=  \exp \left( C p^{\beta} \right), \  C, \ \beta = \const > 0. \eqno(4.13)
$$
 \ Obviously, the r.v. $  \tau $ for which

$$
\forall p \ge 1 \ \Rightarrow \ |\tau|_p \ge \psi_{\beta,C}(p)
$$
does not satisfy the Cramer's condition. \par

\ Let $ \   \xi \ $ be a r.v. belongs to the $ \  G \psi_{\beta,C}(\cdot) \  $ space:

$$
||\xi||  G \psi_{\beta,C} = 1, \eqno(4.14a)
$$
or equally

$$
|\xi|_p \le \exp \left\{ C p^{\beta} \ \right\}, \  p \in [1, \infty). \eqno(4.14b)
$$

 \  The last restriction is quite equivalent to the following tail estimate

$$
T_{\xi}(y) \le \exp \left(  \ - C_1(C, \beta) \ [  \ln(1 + y)   ]^{1 +1/\beta}  \  \right),  \ y > 0. \eqno(4.15)
$$

\vspace{4mm}

\section{Main results: connection between tail behavior and Orlicz's space norm.}

 \vspace{4mm}

 \ We retain the notations and definitions of the previous sections, in particular,

$$
 G(u) = G[\psi](u) =  h^*[\psi](\ln u), \ \psi \in G\Psi  \eqno(5.0)
$$
etc.  Define also the following Young-Orlicz function $ \  N[\psi](u) := $

$$
 \exp[ G(u) ] = \exp \left[ \ h^*[\psi](\ln |u|) \ \right], \ u \ge e; \ N[\psi](u) = C \ u^2, \ |u| < e. \eqno(5.1)
$$

 \ We will prove in this section that the tail estimate (2.2)  of the r.v. $ \ \xi  \ $
 is completely equivalent under some simple conditions to the  finiteness of its Orlicz's norm
 $ \  || \xi|| LN[\psi].  \  $ \par
 \ Recall that we do not suppose that the r.v. $ \ \xi  \ $ satisfies the Cramer's condition. \par

\vspace{4mm}

{\bf Proposition 5.1.}  If  for some  r.v.  $ \ \xi \   $ there holds    $ \ || \xi|| LN[\psi] = K \in (0,\infty), \ $ then

$$
 T_{\xi} ( y) \le \exp \left( \ - h_{\psi}^* (\ln (y/(C \ K)) ) \ \right), \ y \ge e \cdot ||\zeta||, \eqno(5.2)
$$
 \ {\bf Proof} basing only on the Tchebychev-Markov inequality  is at the same as before in the inequality (2.2), see [2], chapters 2,3;
[6], [7],  [19].  Namely,  we  deduce that for some positive finite constant  $ \ C_1 \ $

$$
{\bf E} \exp \left(  G( |\xi|/C_1  \right) < \infty.
$$

 \ It  remains to use the Tchebychev-Markov inequality. \par

\vspace{4mm}

{\bf Proposition 5.2.}   Assume  in addition to the foregoing conditions on the function $  \ \psi(\cdot) \ $
 that the function $  \  G[\psi](u)  \  $ satisfies the following restriction:
$$
\exists \ \alpha = \const \in (0,1), \ \exists K  = \const > 1, \forall  x \in (0,\infty) \ \Rightarrow
G(x/K) \le \alpha \ G(x).  \eqno(5.3)
$$
 \ If  for some  r.v. $ \ \xi \ $

$$
  T_{\xi} ( y) \le \exp \left( \ - h_{\psi}^* (\ln (y) ) \ \right), \ y \ge e, \eqno(5.4)
$$
then the r.v. $ \ \xi \ $ belongs to the Orlicz space $  \  LN[\psi]: $

$$
||\xi||LN[\psi] \le C(\psi,\alpha,K) \ < \infty. \eqno(5.5)
$$

 \vspace{4mm}

 \ {\bf Proof} is more complicated than one for proposition 5.1.  It used the following auxiliary fact.  \par

\vspace{4mm}

 \ {\bf Lemma 5.1.} Let  a function $ \  g: R_+  \to R_+  \  $ be monotonically increasing,  $  \   T = T_{\xi}(x), \
 S = S_{\eta}(x), \ x \ge 0  \  $ be two tail functions correspondingly for non - negative r.v. $ \  \xi, \ \eta \  $ and
such that

$$
T_{\xi}(x) \le  S_{\eta}(x), \ x \ge 0.
$$
 \ We assert:

$$
\int_0^{\infty} g(x) \left|dT_{\xi}(x) \right| \le \int_0^{\infty} g(x) \left|dS_{\eta}(x) \right|.
$$

\vspace{4mm}

{\bf Proof} of lemma 5.1. One can suppose without loss of generality that both the tail functions $ \ T \ $ and $ \ S \ $ are continuous
and strictly  decreasing. Further,  one can realize both the r.v. $  \ \xi, \ \eta \ $  on the classical probability space
$  \ \Omega = \{ \omega\} = [0,1] \ $ equipped with ordinary Lebesgue measure:

$$
\xi = \xi(\omega) = (1 - T)^{-1}(\omega), \ \eta = \eta(\omega) = (1 - S)^{-1}(\omega),
$$
where $ \  f^{-1} \  $ denotes the inverse function. \par
 \ We have $ \  \xi(\omega) \le \eta(\omega)  \ $ a.e., therefore $  \ g(\xi) \le g(\eta) $ a.e., and all the more so

$$
{\bf E} g(\xi) = - \int_0^{\infty} g(x) d T_{\xi}(x) \le   - \int_0^{\infty} g(x) d S_{\xi}(x) = {\bf E} g(\eta),
$$
Q.E.D. \par

  \ {\bf Proof of proposition 5.2.}   Let the pair of numbers $  \  (\alpha, K)  \  $ be from the  condition (5.3).
We have  relaying the proposition of Lemma 5.1 $ \ {\bf E} \exp \left( G(\xi/K) \right) = $

$$
 \int_0^{\infty}  \  \exp G(x/K)  \ |  d  T_{\xi}(x) \ |
 \le \int_0^{\infty}  \ \exp G(x/K)  \ |  d \exp( - G(x))| \le \eqno(5.6)
$$

$$
 \int_0^{\infty}  \ \exp [ \alpha G(x)]  \ |  d \exp( - G(x))| =\int_0^1 z^{-\alpha} d z = \frac{1}{1 - \alpha} < \infty.  \eqno(5.7)
$$
 \ We used by passing $ \ (5.6)  \to (5.7) \ $ the fact quite thin from an article  [14].
 It follows immediately from this estimates that  $ \   \xi \in L(N[\psi]),  \  $ see for example [19],  p. 31 - 33. \par

 \vspace{4mm}

  \  {\bf Examples.} The condition (5.3) is  satisfied  for example for the functions of the form

$$
\psi(p) = p^{1/m} L(p), \ \psi(p) = c p \ [\ln ( p + 1)]^r \ L(\ln p),
$$
where  $ \ m, c = \const \in (0, \infty), \ r = \const \in R \ $ and $  \  L(\cdot) \  $  is positive continuous slowly varying at
infinity function. \par

\vspace{4mm}

 \  {\bf Counterexample.}  The function
$$
\psi_{(r)}(p) = 1, \ 1 \le p \le r;  \ \psi_{(r)}(p) = \infty, \ p > r; \  r = \const > 1,
$$
for which the correspondent function has a form

$$
h^*(\ln u) = r \ln u
$$
does not satisfy the condition (5.3). Actually,  for the r.v.  $  \   \eta \  $  from the space $ \  L_r = L_r(\Omega) \ $
the correspondent tail estimate has a form

$$
T_{\eta} (u) \le c u^{-r},
$$
but the inverse conclusion is not true. \par

\vspace{4mm}

\section{Concluding remarks.}

\vspace{4mm}

 \  {\bf A.} \ It is interest by our opinion to obtain the generalization of  results of this report into multidimensional
case, i.e. into random vectors, alike in the article [7]. \par

\vspace{4mm}

 \ {\bf B. }  We mention  even briefly an important  possible application of obtained results: a Central Limit Theorem in Banach
spaces, in the spirit of [8], section 4.1. \par

\vspace{4mm}

{\it  C.} The case of finite support $ \ \psi: b < \infty. $ \par

 \  In this case approvals 4.1,  5.1, and  5.2 are in general case incorrect. The correspondent counterexamples may be found in the
article [7].  Thus, the problem of description of correspondence  between tail behavior and Grand Lebesgue Space norm  is in this
case an open problem. \par

 \vspace{6mm}

 {\bf References.}

 \vspace{4mm}

{\bf 1. Bennet C., Sharpley R.}  {\it  Interpolation of operators.} Orlando, Academic
Press Inc., (1988). \\

 \vspace{3mm}

{\bf 2.  Buldygin V.V., Kozachenko Yu.V. }  {\it Metric Characterization of Random
Variables and Random Processes.} 1998, Translations of Mathematics Monograph, AMS, v.188. \\

 \vspace{3mm}

 {\bf 3. A. Fiorenza.}   {\it Duality and reflexivity in grand Lebesgue spaces. } Collect. Math.
{\bf 51,}  (2000), 131-148. \\

 \vspace{3mm}

{\bf  4. A. Fiorenza and G.E. Karadzhov.} {\it Grand and small Lebesgue spaces and
their analogs.} Consiglio Nationale Delle Ricerche, Instituto per le Applicazioni
del Calcoto Mauro Picone”, Sezione di Napoli, Rapporto tecnico 272/03, (2005).\\

 \vspace{3mm}

{\bf 5.  T. Iwaniec and C. Sbordone.} {\it On the integrability of the Jacobian under minimal
hypotheses. } Arch. Rat.Mech. Anal., 119, (1992), 129-143. \\

 \vspace{3mm}

{\bf 6. Kozachenko Yu. V., Ostrovsky E.I. }  (1985). {\it The Banach Spaces of random Variables of subgaussian Type. } Theory of Probab.
and Math. Stat. (in Russian). Kiev, KSU, 32, 43-57. \\

\vspace{3mm}

{\bf 7. Kozachenko Yu.V., Ostrovsky E., Sirota L.}  {\it Relations between exponential tails, moments and
moment generating functions for random variables and vectors.} \\
arXiv:1701.01901v1 [math.FA] 8 Jan 2017 \\

  \vspace{3mm}

{\bf 8. Ostrovsky E.I. } (1999). {\it Exponential estimations for Random Fields and its
applications,} (in Russian). Moscow-Obninsk, OINPE. \\

 \vspace{3mm}

{\bf 9. Ostrovsky E. and Sirota L.} {\it Vector rearrangement invariant Banach spaces
of random variables with exponential decreasing tails of distributions.} \\
 arXiv:1510.04182v1 [math.PR] 14 Oct 2015 \\

 \vspace{3mm}

{\bf 10. Ostrovsky E. and Sirota L.}  {\it Non-asymptotical sharp exponential estimates
for maximum distribution of discontinuous random fields. } \\
 arXiv:1510.08945v1 [math.PR] 30 Oct 2015 \\

 \vspace{3mm}

{\bf 11. Ostrovsky E.I.}  {\it About supports of probability measures in separable Banach
spaces.} Soviet Math., Doklady, (1980), V. 255, $ \ N^0 \ $ 6, p. 836-838, (in Russian).\\

 \vspace{3mm}

 {\bf 12. Ostrovsky E. and Sirota L.} {\it Criterion for convergence almost everywhere, with applications.} \\
arXiv:1507.04020v1 [math.FA] 14 Jul 2015. \\

 \vspace{3mm}

{\bf 13. Ostrovsky E. and Sirota L.}  {\it Schl\"omilch and Bell series for Bessel's functions, with probabilistic applications.} \\
 arXiv:0804.0089v1 [math.CV] 1 Apr 2008 \\

 \vspace{3mm}

{\bf 14. Ostrovsky E. and Sirota L. } {\it Sharp moment estimates for polynomial martingales. } \\
arXiv:1410.0739v1 [math.PR] 3 Oct 2014 \\

 \vspace{3mm}

{\bf 15. Ostrovsky E., Rogover E. } {\it Exact exponential bounds for the random field
maximum distribution via the majorizing measures (generic chaining).} \\
 arXiv:0802.0349v1 [math.PR] 4 Feb 2008 \\

 \vspace{3mm}

{\bf 16. Ostrovsky E. and Sirota L. } {\it   Entropy and Grand Lebesgue Spaces approach for the problem  of
Prokhorov - Skorokhod continuity of discontinuous random fields. }\\
arXiv:1512.01909v1 [math.Pr] 7 Dec 2015 \\

 \vspace{3mm}

{\bf 17. Ostrovsky E. and Sirota L. } {\it  Fundamental function for Grand Lebesgue Spaces.  }
arXiv:1509.03644v1  [math.FA]  11 Sep 2015 \\

\vspace{3mm}

{\bf  18. Eugene Seneta.} {\it Regularly Varying Functions.} Lectures Notes in Mathematics, {\bf 508}, (1976). \\

\vspace{3mm}

{\bf 19.  O.I.Vasalik, Yu.V.Kozachenko, R.E.Yamnenko.} {\it $ \ \phi \ - $ subgaussian  random processes. } Monograph, Kiev, KSU,
2008;  (in Ukrainian). \\

\vspace{3mm}

\end{document}